\documentclass[reqno]{amsart}

\usepackage{amssymb, amsmath}

\newtheorem{thm}{Theorem}[section]
\newtheorem{cor}[thm]{Corollary}
\newtheorem{lem}[thm]{Lemma}
\newtheorem{defn}[thm]{Definition} 
\newtheorem{prop}[thm]{Proposition}

\newcommand\ff{{\mathbb F}}
\newcommand\qq{{\mathbb Q}}
\newcommand\zz{{\mathbb Z}}

\newcommand\nc\newcommand
\nc\cO{{\mathcal O}}
\nc\cL{{\mathcal L}}
\nc\ep{\epsilon}
\nc\la{\lambda}
\nc\ze{\zeta}
\nc\om{\omega}
\nc{\al}{\alpha}
\nc{\La}{\Lambda}
\nc{\dbar}{\leftrightarrow}
\nc{\si}{\sigma}
\nc\ab{\stackrel}
\nc\PI{P_{\infty}}
\nc\id{\text{Id}}
\nc\dv{\text{div}}
\nc\phione{\phi^{(1)}}
\nc\phitwo{\phi^{(2)}}
\nc\du{K[\ep]}
\nc\dup{\ff_p[\ep]}
\nc\duq{\ff_q[\ep]}

\begin{document}

\title{A Weil pairing on the $p$-torsion of ordinary elliptic curves over $K[\ep]$}
\author{Juliana V. Belding}
\address{University of Maryland, College Park}
\date{March 29, 2007}
\thanks{This work was made possible with the support of a VIGRE fellowship.}
\thanks{Thanks to Professor Lawrence Washington for his helpful discussions about this work.}

\maketitle

\begin{abstract}
For an elliptic curve $E$ over any field $K$, the Weil pairing $e_n$ is a bilinear map on $n$-torsion. For $K$ of characteristic $p>0$, the map $e_n$ is degenerate if and only if $n$ is divisible by $p$. In this paper, we consider $E$ over the dual numbers $\du$ and define a non-degenerate ``Weil pairing on $p$-torsion" which shares many of the same properties of the Weil pairing.  We also show that the discrete logarithm attacks on $p$-torsion subgroups of Semaev and R\"uck may be viewed as Weil-pairing-based attacks, just like the MOV attack. Finally, we describe an attack on the discrete logarithm problem on anomalous curves, analogous to that of Smart, using a lift of $E$ over $\dup$. 
\end{abstract}

\section{Introduction}

Let $E$ be an ordinary elliptic curve over $K$,  an algebraically closed field of characteristic $p > 0$. For $n$ relatively prime to $p$, the Weil pairing is a bilinear, non-degenerate map
$$e_n: E[n] \times E[n] \rightarrow \mu_n(K)$$ 
 where $E[n] \simeq \zz/n\zz \times \zz/n\zz$ is the $n$-torsion subgroup of $E$ and $\mu_n(K)$ is the group of $n^{th}$ roots of unity of $K$. The Weil pairing is a useful tool in both the theory and application of elliptic curves.  
  
For $p | n$, however, the Weil pairing is degenerate. This is true for two reasons:
$K$ contains no non-trivial $p^{th}$ roots of unity and $E[p] \simeq \zz/p\zz$.  Each of these facts implies that  $e_p(P,Q) = 1$ for all $P,Q \in E[p]$. (The second implies degeneracy since the Weil pairing satisfies the property that $e_n(P,P) = 1$.) 

In this paper, we remedy this situation by considering $E$ over the ring of dual numbers $\du$. Through this deformation of $K$, we find substitutes for the ``missing" geometric points and therefore are able to define a non-degenerate ``Weil pairing" for $n = p$.  In the process, we demonstrate that the discrete logarithm attacks on $p$-torsion subgroups of \cite{ruck} and \cite{sem} are essentially Weil-pairing-based attacks, no different than the MOV attacks on $n$-torsion subgroups for $(n,p)=1$.

In section \ref{prelim}, we give an introduction to elliptic curves over the dual numbers. In sections \ref{milleralg} and \ref{semalg}, we recall Miller's algorithm for computing the Weil pairing and Semaev's algorithm for solving the discrete log problem (DLP) on $p$-subgroups of elliptic curves. In sections \ref{defn} and \ref{proof}, we define the ``Weil pairing on $p$-torsion" $e_p$ over the dual numbers, show its direct relation to Semaev's algorithm, and prove that it satisfies the basic properties of the Weil pairing.  We also describe how $e_p$ can be used to solve the DLP on $p$-torsion subgroups of an elliptic curve. In section \ref{ruckalg}, we give a simple way to compute the pairing using the algorithm of R\"uck defined in  \cite{ruck}. In section \ref{isog}, we describe how the map $e_p$ behaves with respect to isogenies of elliptic curves. In the last section, we give another application of elliptic curves over the dual numbers, namely  a DLP attack on anomalous curves, analogous to that of Smart in \cite{smart}.

\section{Preliminaries}

\subsection{Elliptic Curves over the Dual Numbers}\label{prelim}

The {\bf ring of dual numbers} of the ring $R$ is $R[x]/(x^2)$, denoted $R[\ep]$ with $\ep^2 = 0$. Considering elliptic curves over the dual numbers was proposed in \cite{virat}, where Virat introduced a cryptosystem based on elliptic curves over $\duq$, the dual numbers of $\ff_q$. 

Let $K$ be an algebraically closed field of characteristic $p \neq 0,2,3$.
Let $E$ be the elliptic curve over $K$ given by the Weierstrass equation $y^2 = x^3 + Ax + B$. 
Let $\tilde{A} = A + A_1\ep$ and $\tilde{B} = B + B_1\ep$, for some $A_1,B_1 \in K$.
We call the curve $y^2 = x^3 + \tilde{A}x + \tilde{B}$ a {\bf lift of $E$ to $\du$}, and denote it as  $\tilde{E}$.

The set of points $\tilde{E}(\du)$ consists of two sets: 
\begin{itemize}
\item {\it Affine Points: } $P = (x_0 + x_1\ep: y_0 +y_1\ep:1)$ such that 
\begin{equation}
(x_0,y_0) \in E(K) \text{ and } (2y_0)y_1 = (3x_0^2 + A)x_1 + A_1x_0 + B_1.
\end{equation} 
\item {\it Points at Infinity: } $\cO_k = (k\ep:1:0)$ for all $k \in K$.  

\end{itemize}
Let $\Theta$ denote the set $\{\cO_k | k \in K \}$ and let $\PI$ denote $\cO_0$. The standard addition law for elliptic curves may be extended to give an addition law on $\tilde{E}(\du)$ (see \cite{lcw}, p. 61). An easy calculation shows that 
\[  \begin{matrix}
K^+ & \rightarrow & \Theta \\
k & \mapsto & \cO_k 
\end{matrix}  \]
is an isomorphism. Thus,  $\tilde{E}(\du)$ contains the $p$-torsion subgroup $\Theta$, and there is an exact sequence
$$0 \rightarrow \Theta \rightarrow \tilde{E}(\du) \rightarrow E(K) \rightarrow 0.$$

If $\tilde{A} = A$ and  $\tilde{B} = B$, we call $\tilde{E}$  the {\bf canonical lift} of $E$, since 
the $p$-torsion points $E[p]$ remain $p$-torsion points in $\tilde{E}$. (This terminology comes from the definition of the canonical lift of an elliptic curve to $\qq_q$.)  
For the remainder of the paper (except in Section \ref{smart}), we will assume we are in this situation. In this case, the sequence splits and every point of $\tilde{E}$ may be decomposed as a point of $E(K)$ and a point of infinity. A straightforward calculation using the addition laws gives the following lemma. (Note that  $3x_0^2 + A \neq 0$ for points of order 2, since the curve is non-singular.)

\begin{lem}\label{decomp} Let $\tilde{P} \in \tilde{E}(\du)$ with  $\tilde{P} = (x_0 + x_1\ep: y_0 +y_1\ep:1)$. Then there exists a unique $k \in K$ such that  $\tilde{P} = P + \cO_k$, with $P = (x_0:y_0:1) \in E(K)$. Furthermore  
\[
k =
\begin{cases}
 -\frac{x_1}{2y_0} & \text{if } y_0 \neq 0 \\
 -\frac{y_1}{3x_0^2 + A} & \text{if } y_0 = 0.
\end{cases}
\]
\end{lem}

Note that if $y_0 \neq 0$, the point $(x_1,y_1)$ lies on the line through the origin with slope $\frac{3x_0^2 + A}{2y_0}$, which is precisely the tangent space of the elliptic curve point $(x_0, y_0)$. Thus points of $\tilde{E}(\du)$ may be thought of as points of $E(K)$ with extra ``derivative" information.  (In fact, the set of points of $\tilde{E}(\du)$ may be naturally identified with the tangent bundle of the variety $E$.)

The canonical lift $\tilde{E}$ has $p$-torsion $\tilde{E}[p] = E[p] \oplus\Theta$. Furthermore, $\mu_p(\du)$ has non-trivial $p^{th}$ roots of unity, in particular the subgroup $\{ 1 + a\ep: a \in K \}$. 
Thus we will see that there is a non-degenerate ``Weil pairing" on the $p$-torsion of $\tilde{E}$. Before we proceed we recall Miller's method of computing the Weil pairing. 

\subsection{Miller's algorithm for computing the Weil pairing}\label{milleralg}

Let $(n,p) = 1$. Let $P,Q \in E[n]$, and let $D_P$, $D_Q$ be divisors with disjoint support which sum to $P,Q$ respectively.  Let $f_P, f_Q$ be functions with divisors $nD_P, nD_Q$ respectively. The Weil pairing is defined as 
$$e_n (P,Q) = \frac{f_P(D_Q)}{f_Q(D_P)}.$$ 

This definition is independent of the choices of divisors by Weil reciprocity. In \cite{miller}, Miller gives a way to compute the value $f_P(D_Q)$. As this will be the foundation for the definition of the ``Weil pairing on $p$-torsion," we recall the details here. 

Let $P \in E[n]$. Choose any two points $T, R \in E(K)$ such that the divisors $D_P = (P+T) -(T)$ and $D_Q = (Q + R) - (R)$ are disjoint. Let $f_P$ be the function with divisor $nD_P$. Note that this function is unique only up to a non-zero constant. Following \cite{miller}, in such situations, we choose the unique function with the value $1$ at $\PI$, which we call the {\bf normalized} function. (Note that since we are calculating the ratio $f_P(Q + R)/f_P(R)$, such constants may in fact be disregarded.)

Let  $f_k$ denote the (normalized) function with $$\dv(f_k) = k(P +T) - k(T) - (kP + T) + (T).$$ Note that $\dv(f_1) = 0$, so $f_1 \equiv 1$ . Also note that $\dv(f_P) = \dv(f_n)$  and $\dv(f_{i+j}) = \dv(f_if_jh_{i,j})$ where $$\dv(h_{i,j}) =   - ((i+j)P + T) + (iP + T) + (jP + T) - (T).$$
Thus $f_P(Q) = f_n(Q)$ can be calculated recursively by using an addition chain decomposition for $n$. 

An addition chain  for a positive integer $n$ is an increasing sequence of  integers $S \subset \{ 1, ..., n\}$ such that for each $k \in S$  with $k>1$, there exist $i,j \in S$ such that $i + j= k$. Given an addition chain $S$, an {\bf addition chain decomposition $C$} of $n$ is a sequence of steps of the form $(k \mapsto i,j)$ with $i+j = k$ and $i,j,k \in S$ which decomposes $n$ into the sum of $n$ ones: $\underbrace{1+...+1}_n$. 
Note that any decomposition will consist of exactly $n-1$ steps. 

Thus, since $f_k(Q) = f_i(Q)f_j(Q)h_{i,j}(Q)$ and  $f_1 \equiv 1$, $f_n(Q) $ will be the product of $n-1$ contributions of the form $h_{i,j}(Q)$. For example, if $n = 11$ and $S =\{1,2,4,8, 10, 11\}$, then one possible decomposition is
\[f_{11} = f_1f_{10}h_{1,10} = f_1f_2f_8h_{1,10}h_{2,8}  = ... = f_1^{11}h_{1,10}h_{2,8}h_{4,4}h_{2,2}^2h_{1,1}^5.\]
Given an addition chain decomposition $C$ for $n$, we write $\prod_{C} h_{i,j}(Q)$ to denote the value $f_n(Q)$. Note that there always exists a decomposition with $O(\log n)$ distinct $h_{i,j}$.

Let $\ell_{i,j}$ denote the line through $iP$ and $jP$, and let $v_i$ denote the vertical line through $iP$. Note that 
\[ \dv(\ell_{i,j}) = (iP) + (jP) + (-(i+j)P) - 3\PI \text{   and    } \dv(v_i) = (iP) + (-iP) - 2\PI. \]
Let $\tau$ denote translation by $-T$. Then

\begin{equation}\label{h}
h_{i,j} = 
\begin{cases}
  \frac{\ell_{i,j}}{v_{i+j}}\circ \tau   & \text{$i+j \neq n $}, \\
  v_i \circ \tau    & i+j = n.
\end{cases}
\end{equation}

As is remarked in \cite{freyruck}, this calculation of $f_P(Q)$ may be interpreted as exponentiation  in a generalized jacobian with modulus $(Q+T) - (T)$. The $h_{i,j}$ are simply cocycle values.  A good source for this viewpoint is \cite{dechene}.  

\subsection{Semaev's algorithm for solving the DLP on anomalous elliptic curves}\label{semalg}

Let $K = \ff_q$ be a finite field of characteristic $p$. In \cite{sem}, Semaev proposed a polynomial time algorithm for solving the DLP on elliptic curves over $K$ which contain a point of order $p$, using the following map:
\[
\begin{matrix}
\la: &  E[p] & \rightarrow & K^+  \\ 
      & P       & \mapsto      & \frac{f'_P}{f_P}(R)  \\
      & \PI      & \mapsto     & 0 
\end{matrix}
\]
where $D_P$ is any divisor of degree 0 which sums to $P$, $f_P$ is any function with $\dv (f_P) = pD_P$, and $R \in E[p]$ with $R \neq \PI$. Here $f'_P$ denotes $\frac{d}{dx}f_P$. 

To see how this map is used to solve the DLP, consider $P,Q \in E[p]$ with $Q = nP$. Using the standard $\log p$ addition chain decomposition, we can compute $\la(P), \la(Q)$ in time $O(\log p)$, and then solve $n\la(P) = \la(Q)$ for $n \in K^+$ by Euclid's algorithm.

\begin{prop}\label{semprop} (Semaev, \cite{sem}) The map $\la$ is defined and non-zero for any $R \notin E[2]$. Furthermore, $\la$ is an injective homomorphism with respect to $P$ and is independent of the divisor $D_P$. 
\end{prop}

This is proved explicitly in \cite{sem} and in fact, the proof holds for any algebraically closed field $K$ of characteristic $p > 0$.  The proposition also follows  from considering the map:
\[  
\begin{matrix}
E[p] & \rightarrow & Pic^0_K(E)[p] & \rightarrow & \Omega^h_K(E) & \rightarrow & \cL_{div(dt)} & \rightarrow  & K^+  \\
P & \ab{\delta}\mapsto & D_P & \ab{\rho}\mapsto & \frac{df_P}{f_P} & \ab{\psi}\mapsto & \frac{df_P}{dtf_P} & \ab{\varphi}\mapsto & \frac{df_P}{dtf_P}(R) \\ 
\end{matrix}  \] 
where $Pic^0_K(E)$ is the group of divisor classes of $E$ of degree 0, $\Omega^h_{K}(E)$ are the holomorphic differentials of the one-dimensional $K(C)$-vector space of differentials,  $ \cL_{div(dt)}$ is the one-dimensional $K$-vector space of functions $g$ with $\dv(gdt) \geq 0$ and $t$ is a uniformizer for the point $R$.

This is an injective homomorphism since  $\rho$ is an injective homomorphism (see \cite{serre}) and $\delta, \psi, \varphi$ are isomorphisms. This is noted in \cite{ruck}, where the attack on the DLP is extended to the $p$-subgroup of the divisor class group of a curve of arbitrary genus.

The computation method proposed in \cite{sem} is a variation on Miller's algorithm. Let $T$ be a point of order two and let $R \in E[p]$. Let $f_Q$ be the function with $\dv(f_Q) = D_Q = (Q+T) - (T)$. As in section \ref{milleralg}, the value of the function $\la$  may be computed by using an addition chain decomposition and summing contributions of the form  $\frac{h_{i,j}'}{h_{i,j}}(R)$, where $h_{i,j}$ is as in Section \ref{milleralg}. That is, $\frac{f'_P}{f_P}(R) = \sum_{C} \frac{h_{i,j}'}{h_{i,j}}(R)$. (Remark: In \cite{sem}, the function $\frac{v_iv_j}{l_{-i,-j}}$ is used, which is equivalent up to constant since it has the same divisor.) 

To compute $h' = \frac{d}{dx}h$, we make use of the invariant differential property $\frac{dx}{y} \circ \tau = \frac{dx}{y}$. Let $g$ be a function expanded in a power series around $x$. Then $g \circ \tau$ can be expanded in a series around $x \circ \tau$ with the same coefficients, and so
$$\frac{d(g \circ \tau)}{dx} = \frac{d(g \circ \tau)}{d(x \circ \tau)} \frac{d(x \circ \tau)}{dx} = (\frac{dg}{dx} \circ \tau) \frac{y\circ \tau}{y}.$$

\noindent Therefore, for $h = \frac{\ell}{v} \circ \tau$, $h'(R) = \frac{y(R+T)}{y(R)}(\frac{\ell}{v})'(R+T)$. (Since $T$ is order 2, translation by $T$ and $-T$ are the same.) \\

\noindent {\bf Remark: } The choice of the divisor $D_P = (P+T) - (T)$ avoids any possible zeros or undefined values  when evaluating the lines through multiples of $P$ at $R$, which is itself a multiple of $P$. Note that when $p>7$, for a fixed point $P$, it is always possible to choose an $R \in E[p]$ such that the lines in a $\log p$ addition chain decomposition will not have $R$ as a zero. However,  since the homomorphism $\la$ is {\it not} independent of $R$, in order to have it well-defined it is necessary to choose an evaluation point that works for all $P$, which explains Semaev's use of a translation point.

As is the case for Miller's algorithm to compute the Weil pairing, this calculation may be interpreted as exponentiation in a generalized jacobian, after a slight modification. Note that if we use the divisor $D_P = (P) - (\PI)$, the $h_{i,j}$ are simply ratios of lines through multiples of $P$, and thus evaluating at $R+T \notin E[p]$ gives well-defined, non-zero values. In this case, we may calculate the value $\frac{f_P'}{f_P}(R+T)$  using exponentiation in a generalized jacobian with modulus $2(R+T)$ for $R \in E[p]$, with $T$ of order 2.  The value will differ from the value $\la(R)$ by the constant factor $y(R+T)/y(R)$. 

\section{A ``Weil pairing" on the $p$-torsion of $\tilde{E}(\du)$}\label{defn}

Let $\tilde{E}$ denote $\tilde{E}(\du)$, the canonical lift of $E: y^2 = x^3 + Ax + B$ to $\du$. We define the pairing $$e_p: \tilde{E}[p] \times \tilde{E}[p] \rightarrow \mu_p(\du)$$ by first defining a bilinear map $e$ on $E[p] \times \Theta$, and then extending it to $\tilde{E}[p]$ in such a way that the necessary properties hold. 

Let $P \in E[p]$ and let $T$ be a point of order two. Consider the divisor $D_P = (P+T) -(T)$. Let $f_P$ be the function on $E$ with divisor $pD_P$, unique up to a non-zero constant. We use the notation of section \ref{milleralg}. Recall that to compute $f_P$ evaluated at a point $Q$, we choose an addition chain decomposition for $p$ and compute the product of cocycle contributions of the form $h_{i,j}(P)$, where $h_{i,j}$ are ratios of lines translated by $T$.  

Any function in $K(E)$ is a well-defined function on the affine points of $\tilde{E}(\du)$, provided that the denominator is invertible. We will see that this is true for $h_{i,j}$ on certain points of $\tilde{E}$, thereby making the computation of $\prod_{C} h_{i,j}$ legitimate. 

\begin{defn}
Fix $R \in E[p]$ such that $R \notin E[2]$. Let $C$ be an addition chain decomposition for $p$. Define the map $e: E[p] \times  \Theta \rightarrow \mu_p(\du)$ by
$$e(P,\cO_k) = 
\begin{cases}
\prod_{C} \frac{h_{i,j}(\cO_k + R)}{h_{i,j}(R)} & \text{ if } P, \cO_k \neq \PI   \\
1 & \text{ if } P = \PI \text{ or } \cO_k = \PI \\
\end{cases}
$$
\end{defn}

The proof  of the following theorem is given in the next section.

\begin{thm}\label{mainthm} The map $e$ is well-defined and bilinear and independent of the addition chain decomposition of $p$. Furthermore, for any divisor $D_P$ summing to $P$ and any $R \in E(K)$, 
$$
 e(P, \cO_k)  = 1 - 2\Big(y\frac{f'_P}{f_P}\Big)(R)k\ep,
 $$
where $f_P$ is the normalized function with divisor $pD_P$.
\end{thm}

We now may define the {\bf Weil pairing on $p$-torsion.} 
Extend the map $e$ to $$e_p: \tilde{E}[p] \times \tilde{E}[p] \rightarrow \mu_p(\du)$$ such that 
\begin{itemize}
\item $e_p(P, \cO_k) = e(P, \cO_k)$ for all $P \in E[p]$,
\item $e_p(P,Q) =1$, for all $P,Q \in E[p]$, 
\item $e_p(\cO_k,\cO_j) = 1$, for all $j,k \in K$,
\item $e_p$ is bilinear, 
\item $e_p$ is anti-symmetric: $e(P,Q) = e(Q,P)^{-1}.$
\end{itemize}

\begin{thm}\label{nondeg}
The map $e_p$ is non-degenerate. That is, if $e_p(P, Q) = 1$ for all $P \in \tilde{E}[p]$, then $Q = \PI$, and if $e_p(P, Q) = 1$ for all $Q \in \tilde{E}[p]$, then $P = \PI$. 
\end{thm}

The proof of this theorem is given in the next section. \\

\noindent{\bf Remark: } 
Note that we are defining $e_p(P,\cO_k)$ to be the result of Miller's algorithm to compute
$$ \frac{f_P(\cO_k + R)}{f_P(R)}.$$
This definition can thus be viewed as the analog of the Weil pairing definition for $n$ prime to $p$:
$$e_n (P,Q) = \frac{f_P(Q + R)}{f_P(R)}\frac{f_Q(R)}{f_Q(P + R)}.$$ 
Recall that Miller's algorithm computes the value of $f_Q$ as the product of ratios of lines through multiples of the point $Q$. For $Q = \cO_k$, this involves products of  lines through points at infinity (which would then be evaluated at affine points of $E(K)$). 
Assuming such a line has the form $\ell = 0$ with $\ell(x,y,z) = ax + by + cz$ and $a,b,c \in \du$, there is not a unique choice for such a line. For example, any line of the form $a\ep x+ cz$, for $a \in K, c\in \du$, passes through the points $\cO_k, \cO_j$. We make the choice of the line $\ell = cz$. When evaluated at affine points, this becomes the constant function $c$ which normalized is just $1$.  The value of $\frac{f_Q(R)}{f_Q(P + R)}$ for $Q = \cO_k$ may therefore naturally be considered to be 1.  \\

We now show how the Weil pairing $e_p$ can be used to solve the DLP on $p$-subgroups of elliptic curves over $\ff_q$. Given $P,Q \in E[p]$ with $Q = nP$, calculate $e_p(P, \cO_1) = 1 + a\ep$ and $e_p(Q, \cO_1) = 1 + b\ep$. Since $e_p$ is bilinear, $e_p(Q, \cO_1) = e_p(P, \cO_1)^n = (1 + a\ep)^n = 1 + na\ep$. Thus it suffices to solve the equation $b = na$ in $\ff_q^+$ for $n \in \zz/p\zz$ by computing the multiplicative inverse of $a$. By Theorem \ref{mainthm}, for $R \in E[p]$, this process is essentially Semaev's algorithm to solve the DLP in $p$-subgroups. Therefore, we see that Semaev's algorithm may be interpreted as a Weil-pairing based attack.

\section{Proof of properties of the pairing }\label{proof}

To show that $e$ is well-defined and bilinear, we relate its calculation to the map $\la$ from section \ref{semalg}. For this, we need the following lemma. 

\begin{lem}\label{key}
Let $\ell_{i,j}$ denote the line through $iP$ and $jP$, and let $v_i$ denote the vertical line through $iP$. Let $\tau$ denote translation by $-T$ and let 
$$
h_{i,j} = 
\begin{cases}
  \frac{\ell_{i,j}}{v_{i+j}}\circ \tau   & i+j \neq p, \\
  v_i \circ \tau    & i+j = p.
\end{cases}
$$
Let $R \in E[p]$ with $R \neq \PI$. Then
$$\frac{h_{i,j}(\cO_k+R)}{h_{i,j}(R)} = 1 - 2y(R)\frac{h_{i,j}'}{h_{i,j}}(R)k\ep.$$ 
\end{lem}

\noindent {\bf Proof:} We first show that
$$h_{i,j}(\cO_k+R)= h_{i,j}(R) - 2y(R)h_{i,j}'(R)\ep.$$
We can think of this as analogous to the calculus approximation of $f(x_0 + \ep)$ by the value $f(x_0) +  f'(x_0)\ep$.

Let $S = R+T = (x_0, y_0)$. Assume $i+j \neq p$. Fix $i,j$ and let $h_{i,j} = h = \frac{\ell}{v} \circ \tau$. Since we are evaluating $h_{i,j}$ at affine points, we have $\ell = y - mx - b$ and $v = x - c$ for some $m,b,c \in K$. 

Since $v$ is a line through a multiple of $P$, and $S \notin E[p]$, we see that $ x_0 - c \neq 0$. Thus $h(R) = \frac{\ell}{v}(S)$  
is well-defined. Furthermore, since $\cO_k+ S = (x_0 - 2y_0k\ep: y_0 - (3x_0^2 + A)k\ep:1)$, the denominator of $h(\cO_k + R)$ is invertible, and thus the value $h(\cO_k + R)$ is well-defined. Then

\[  \begin{matrix}
h(\cO_k+R) & =  & 
 \frac{\ell}{v}(\cO_k+S) \\
& = & \Big[(y_0 - mx_0- b) + (2y_0m - (3x_0^2 + A))k\ep\Big] \Big/ \Big[(x_0 - c) - 2y_0k\ep\Big] \\
& = & \Big[(y_0 - mx_0- b) + (2y_0m - (3x_0^2 + A))k\ep\Big]\Big[(x_0-c)^{-1} + (x_0-c)^{-2}2y_0k\ep\Big] \\
& = & h(R) + \Big[(2y_0m - (3x_0^2 + A) )(x_0 - c)^{-1} + h(R)(x_0 - c)^{-1}2y_0\Big]k\ep. 
\end{matrix}  \]

\noindent Recall from section \ref{semalg} that $h'(R) = \frac{y(S)}{y(R)} (\frac{\ell}{v})'(S)$. 
Since $v' = 1$ and $l'(S) = \frac{3x_0^2 + A}{2y_0} - m$, we have
$$h'(R) = \frac{1}{2y(R)} \Big( (3x_0^2 + A - 2y_0m)(x_0 - c)^{-1} - 2y_0h(R)(x_0 - c)^{-1} \Big)$$
and therefore  $h(\cO_k+R)  = h(R) - 2y(R)h'(R)k\ep.$ 

For $i+ j = p$, we have $h = v \circ \tau$ and $h'(R) = \frac{y(S)}{y(R)}$ by the equation in Section \ref{semalg}. Then
$$h(\cO_k+R) = v \circ \tau(\cO_k+R) = v(\cO_k+S) = (x_0 - c) - 2y_0k\ep =h(R) - 2y(R)h'(R)k\ep. $$

\noindent It remains to show that $h(R) \neq 0$. The fact that $R \in E[p]$ implies that $S$ is not a zero of the line described by $\ell$ or $v$. Therefore, in both cases, $h(R) \neq 0$, and the result follows. $\Box$ \\

Now we can prove Theorem \ref{mainthm} and \ref{nondeg}.\\

\noindent {\bf Proof (Thm. \ref{mainthm}):} 
Fix $P, \cO_k$ and an addition chain decomposition $C$ for $p$.  Note that by Lemma \ref{key}, $e(P, \cO_k)$ is well-defined. Let $f_P$ be the function with divisor $pD_P$ for $D_P = (P+T) - (T)$. Let $R \in E[p]$ with $R \neq \PI$. Then $\frac{f'_P(R)}{f_P(R)} = \sum_C \frac{h_{i,j}'}{h_{i,j}}(R)$. We have
\[  \begin{matrix}
e(P, \cO_k) & = & \prod_{C} \frac{h_{i,j}(\cO_k+R)}{h_{i,j}(R)} \\
& = & \prod_{C} \Big( 1 - 2y(R)\frac{h_{i,j}'}{h_{i,j}}(R)k\ep \Big) \\
& = & 1 - 2y(R) \Big( \sum_{C} \frac{h_{i,j}'}{h_{i,j}}(R) \Big) k\ep \\
& = & 1 - 2y(R)\frac{f'_P(R)}{f_P(R)}k\ep.
\end{matrix}  \]
Note that $\frac{f'_P(R)}{f_P(R)} = \la(R)$, where $\la$ is the homomorphism with respect to $P$  from section \ref{semalg}. Thus since $e(P, \cO_k) = 1 - 2k(y\frac{f'_P}{f_P})(R)\ep$, the map $e$ is linear in the first coordinate. Furthermore, since $\cO_k + \cO_j = \cO_{k+j}$, we have that $e$ is linear in the second coordinate. Therefore, $e$ is bilinear.
Since $\frac{f_P'(R)}{f_P(R)}$ is independent of addition chain decomposition, so is the value of $e$.  

As shown in \cite{sem}, $\dv(\frac{f_P'}{f_P}) = \dv(\frac{1}{y})$, thus $y \frac{f_P'}{f_P}$ is a constant function on $E(K)$. 
For $R \in E[p]$, we've just seen that $e(P, \cO_k) = 1 - 2(y\frac{f_P'}{f_P})(R)k\ep$. Therefore $$e(P, \cO_k)  = 1 - 2k\Big(y\frac{f'_P}{f_P}\Big)(R)\ep$$
for all $R$, and thus $e$ is independent of $R$. Furthermore, since $\frac{f_P'}{f_P}$ is independent of the divisor for $P$, as shown in \cite{sem}, the value of $e$ is independent of choice of the divisor for $P$.   $\Box$ \\

\noindent {\bf Proof (Thm. \ref{nondeg}):} Let $P \in \tilde{E}(\du)[p]$. We show that if $P \neq \PI$, then there exists $Q \in \tilde{E}(\du)[p]$ such that $e_p(P,Q) \neq 1$. This shows non-degeneracy in the first coordinate, and by the property of anti-symmetry, non-degeneracy in the second coordinate will follow.  

By Lemma \ref{decomp}, $P$ may be written as $P_0 + \cO_k$ for $P_0 \in E(K)$ and $k \in K$. 
If $P_0 \neq \PI$, let $Q = \cO_1$. Then $e_p(P,Q) = e_p(P_0, \cO_1)e_p(\cO_k,\cO_1) = e_p(P_0, \cO_1)$. Let $R \in E[p]$ with $R \neq \PI$. By Proposition \ref{semprop}, $\frac{f'_P}{f_P}(R)$ is non-zero.  Therefore, since $e_p(P_0, \cO_1) = 1 - 2(y\frac{f'_P}{f_P})(R)\ep$ and $R \notin E[2]$, we have that $e_p(P,Q) \neq 1$. 

If $P_0 = \PI$, then $k \neq 0$, since $P \neq \PI$. Let $Q, R \in E[p]$ with $Q, R \neq \PI$. 
Then $e_p(P,Q) = e_p(\cO_k,Q) = 1 + 2(y\frac{f'_Q}{f_Q})(R)k\ep$. Since $k \neq 0$ and $R \notin E[2]$, we have that $e_p(P,Q) \neq 1$, as desired. $\Box$ 

\section{R\"uck's algorithm for solving the DLP on $p$-torsion}\label{ruckalg}

Recall the homomorphism from Section \ref{semalg}:
\[  \begin{matrix}
E[p] & \rightarrow & Pic^0_K(E)[p] & \rightarrow & \Omega^h_K(E) & \rightarrow & \cL_{div(dt)} & \rightarrow  & K^+ \\
P & \ab{\delta}\mapsto & D_P & \ab{\rho}\mapsto & \frac{df_P}{f_P} & \ab{\psi}\mapsto & \frac{df_P}{dtf_P} & \ab{\varphi}\mapsto & \frac{df_P}{dtf_P}(R) \\ 
\end{matrix}  \]

Choosing the  divisor $D_P = (P) - (\PI)$ and evaluation point $R = \PI$, we may compute the value of $\frac{df_P/dt}{f_P}(R)$  by simply summing the slopes of lines through multiples of $P$ for any addition chain decomposition. This fact is noted in \cite{kk}, where it is referred to as the ``R\"uck algorithm,"  and a slight variation is found in \cite{zp}. In \cite{ruck}, R\"uck refers to the result of this algorithm as ``the additive version of the Tate pairing." We make this remark explicit by relating the algorithm to the pairing of $E[p]$ and $\Theta$ which we've defined. 

\begin{prop}\label{sumofslopes} Let $m_{i,j}$ be the slope of the line through $iP$ and $jP$, and let $C$ be an addition chain decomposition for $p$. Then
$$e(P, \cO_k) = 1 + \Big[ \sum_{C} m_{i,j} \Big] k\ep.$$ 
\end{prop}

\noindent{\bf Proof:} As $e$ is independent of divisor and evaluation point, we may choose the divisor $D_P = (P) - (\PI)$ and evaluation point $R = \PI$. This means we must calculate
$$e(P, \cO_k) = 1 - 2\Big(y\frac{f'_P}{f_P}\Big)(\PI)k\ep.$$

Since we evaluate at $\PI$, we want to expand functions around the uniformizer for $\PI$, namely $t  = -\frac{x}{y}$. Using the fact that $\frac{dt}{dx} = \frac{x^3 + Ax + 2B}{2y^3}$,  we are looking to compute
$$\frac{df_P/dt}{f_P}\frac{x^3 + Ax + 2B}{y^2}(\PI).$$ 
Recall that $x$ and $y$ have poles at $\PI$ of order 2 and 3, respectively. In particular, $x = \frac{1}{t^2} + O(t)$ and $y = -\frac{1}{t^3} + O(t)$ (\cite{sil}, p. 113). Thus  $\frac{x^3 + Ax + 2B}{y^2} = -1 + O(t)$, and hence this contributes a factor of $-1$ when we evaluate at $\PI$. 

We now focus on computing $\frac{df_P/dt}{f_P}(\PI)$. Since $D_P = (P) - (\PI)$, this reduces to computing $\frac{dh_{i,j}/dt}{h_{i,j}}$ where $h_{i,j}$ is defined as in  section \ref{semalg}. In particular, we show that
\begin{equation}\label{slope}
\frac{dh_{i,j}/dt}{h_{i,j}} = 
\begin{cases}
   -\frac{1}{t} - m_{i,j} + O(t)    & \text{ if } i + j \neq p, \\
   -\frac{2}{t} + O(t)     & \text{ if } i + j = p 
\end{cases}
\end{equation}

For $i + j \neq p$, $$h_{i,j}  =\frac{\ell}{v} =  \frac{y - mx - b}{x-c} = 
\frac{1}{t} - m + O(t).$$
Thus 
$$\frac{dh/dt}{h} = 
-\frac{1}{t} - m + O(t).$$
For $i+j = p$, $h_{i,j} = v$. Expanding $v$ around $t$, we get $v = x - c = \frac{1}{t^2} - c + O(t)$. Thus $$\frac{dv/dt}{v} = 
-\frac{2}{t} + O(t),$$ 
and (\ref{slope}) is proved. 
 
Note that using an addition chain decomposition for $p$ to calculate $f_P$ will result in $(p-1)$ terms of the form $h_{i,j}$ with exactly one such that $i+j = p$. Thus the pole contributions of the $h_{i,j}$ total to zero in characteristic $p$ and 
$$\frac{df_P/dt}{f_P} = -\frac{p}{t} - \sum_{C} m_{i,j} + O(t) = - \sum_{C} m_{i,j} + O(t) $$
Evaluating at $\PI$ yields the result. $\Box$ 

\begin{cor}
(R\"uck, \cite{ruck}) Let $f_P$ be any function with divisor $p(P) - p(\PI)$ and let $t = -\frac{x}{y}$. Let $m_{i,j}$ denote the slope of the line through $iP$ and $jP$, and let $C$ be an addition chain decomposition for $p$. Then
$$\frac{df_P/dt}{f_P}(\PI) = -\sum_{C} m_{i,j}.$$
\end{cor}

\section{The map $e_p$ and isogenies of $\tilde{E}(\du)$}\label{isog}

Let $\phi:E_1 \rightarrow E_2$ be an isogeny between curves given by the Weierstrass form $y^2 = x^3 + A_ix + B_i$. Let $\tilde{E}_i$ denote the canonical lift of $E_i$, as defined in \ref{prelim}. In this section, we show how to extend $\phi$ to a homomorphism $\tilde{\phi}: \tilde{E}_1 \rightarrow \tilde{E}_2$ in such a way that the following proposition holds:
\begin{prop}\label{isogprop}
For any isogeny $\phi:E_1 \rightarrow E_2$, $$e_p(\tilde{\phi}(P), \tilde{\phi}(Q)) = e_p(P,Q)^{\deg \phi}.$$\end{prop}

As $\tilde{E}_i \simeq E_i \oplus \Theta_i$, it suffices to define $\phi: \Theta_1 \rightarrow \Theta_2$ and then extend it linearly to a map $\tilde{\phi}:\tilde{E}_1 \rightarrow \tilde{E}_2$. 
Let $x_i,y_i$ denote the coordinate functions of $E_i$, and let $t_i = -\frac{x_i}{y_i}$ be a uniformizer at $\PI^{(i)}$, the point at infinity of $E_i$. Let $m \in K$ be such that $t_2 \circ \phi = mt_1 + O(t_1^2)$. 
(To obtain the value $m$, expand $x_1$ and $y_1$  around  $t_1$ and use the fact that  $x_2$ and $y_2$ are rational functions of $x_1$ and $y_1$ to obtain $t_2 \circ \phi$ as a function of $t_1$.)

\begin{defn}\label{isogdef}  For $\phi:E_1 \rightarrow E_2$, define $\tilde{\phi}:\Theta_1 \rightarrow \Theta_2$ by $\tilde{\phi}(\cO_k) = \cO_{mk}.$ 
\end{defn} 

First note that $\tilde{\phi}: \Theta_1 \rightarrow \Theta_2$ is a homomorphism with respect to this definition. 
Furthermore,  it is compatible with composition of isogenies. That is, if $\phi: E_1 \rightarrow E_2$ and $\psi:E_2 \rightarrow E_3,$ are isogenies, then $(\tilde{\psi} \circ \tilde{\phi} )(\cO_k) = \tilde{\psi} (\tilde{\phi} (\cO_k))$. 
This follows from the fact that if $t_2 \circ \phi = m_1t_1 + O(t_1^2)$ and $t_3 \circ \psi = m_2t_2 + O(t_2^2)$, then $t_3 \circ (\psi \circ \phi) = (m_1m_2)t_1 + O(t_1^2)$.

The motivation for the definition is as follows. If $\phi$ inseparable, then $\phi = \phi_s \circ \pi^r$, where $\phi_s$ is separable and the degree of inseparability of $\phi$ is $p^r$. The map $\pi:(x:y:z) \mapsto (x^p:y^p:z^p)$ is well-defined on the points of $\tilde{E}(\du)$, and clearly $(k\ep:1:0) \ab{\pi}\mapsto (0:1:0)$. Thus we should define $\phi(\cO_k) = \PI$. (Note that this agrees with the idea that $\Theta$ is acting as the replacement for the ``missing" geometric points of $p$-torsion, the ``kernel of Frobenius." )  But $m = 0$ if $\phi$ is inseparable, since the order of $t_2 \circ \phi$ at $\PI^{(1)}$ is 
equal to the degree of inseparability (\cite{sil}, p. 76), so $\phi(\cO_k) = \cO_{mk} = \PI^{(2)}$. 

Now consider $\phi$ separable. Then $t_2 \circ \phi$ is a uniformizer for $\PI^{(1)}$, so $m \neq 0$.  
Suppose we want $\phi(\cO_k) = \cO_j$, for some $j \in K$. Since $t_2 \circ \tilde{\phi}(\cO_k) = t_2((j\ep:1:0)) = -j\ep$ and $(mt_1 + O(t_1^2))(\cO_k) = -mk\ep$, it makes sense to define $j = mk$. 

Next we  extend the isogeny $\phi: E_1 \rightarrow E_2$ to the affine points of $\tilde{E}(\du)$. 

\begin{defn}\label{affine}
Let $\tilde{P}$ be a lift of an affine point $P \in E_1(K)$. Let $T \in E_1(K)$ with $T \notin \ker \phi$. Let $\tau$ denote translation and let $\phi_T = \tau_{\phi(-T)} \circ \phi \circ \tau_T$. Define
\[
\tilde{\phi}(\tilde{P}) = 
\begin{cases}
\phi(\tilde{P}) & \text{ if } P \notin \ker \phi \\
 \phi_T(\tilde{P}) & \text{ if } P \in \ker \phi 
\end{cases}
\]
\end{defn}

Note that this definition is independent of $T$. That is, for $T, T' \notin \ker \phi$, we have $\phi_T(\tilde{P}) = \phi_{T'}(\tilde{P})$ for all $P \in \ker \phi$. Furthermore, $\phi_T(\tilde{P}) = \phi(\tilde{P})$ for all $\tilde{P} \in \tilde{E_1}(\du)$ for which both isogenies are defined. This follows from the fact that $\tilde{\phi}$ is a homomorphism (see  Proposition \ref{hom} below).

We need to establish that this definition yields well-defined points in $\tilde{E}_2(\du)$.
First note that any isogeny $\phi$ can be written as $\phi(x,y) = (r(x), ys(x)))$ where $r,s$ are rational functions of $x$ (\cite{lcw}, p. 47). Thus, we can evaluate $\phi$ on affine points $\tilde{P} = (x_0 + x_1k\ep: y_0 +y_1k\ep:1) \in \tilde{E}(\du)$ provided that the denominators of $r$ and $s$ are invertible when evaluated at $\tilde{P}$.  This will be the case for all $\tilde{P}$ such that $P = (x_0, y_0)$ is not a kernel point of $\phi$.

Note that $\phi(\tilde{P})$ is in fact a point of $\tilde{E_2}(\du)$. For all $x,y \in E_1(K)$, 
\begin{equation}\label{image1}
y^2s(x)^2 = r(x)^3 + A_2r(x) + B_2.
\end{equation}
Therefore, for all $x \in K$, we have
\begin{equation}\label{image2}
(x^3 + A_1x + B_1)s(x)^2 = r(x)^3 + A_2r(x) + B_2.
\end{equation}
This is an identity in the function field $K(x)$.   Let $\tilde{P} = (x,y) \in E_1(\du)$. Then $x$ satisfies (\ref{image2}) and since $y^2 =x^3 + A_1x + B_1$, we have that $(x,y)$ satisfies (\ref{image1}).
Therefore $\phi(\tilde{P}) \in \tilde{E_2}(\du)$.

Now consider $\tilde{P}$ such that $P \in \ker \phi$. If $P \neq \PI$ and $\tilde{P} \neq P$, then $\tau_{\phi(-T)} \circ \phi \circ \tau_T$ is well-defined on these $\tilde{P}$ and yields a point of $\tilde{E}_2(\du)$, since translation by $T$ is a map of the curve to itself. 

Combining definitions \ref{isogdef} and \ref{affine}, we can extend $\phi$ to the map $\tilde{\phi}:\tilde{E}_1 \rightarrow \tilde{E}_2$, which we show is a homomorphism. 

\begin{prop}\label{hom}
Let $P,Q \in \tilde{E_1}(\du)$. Then 
$$\tilde{\phi}(P) + \tilde{\phi}(Q) = \tilde{\phi}(P + Q).$$
\end{prop}

\noindent{\bf Proof:} 
By Lemma \ref{decomp}, any point of $\tilde{E_1}(\du)$ decomposes as a point of $E_1(K)$ and $\Theta$. Thus, since $\tilde{\phi}$ is homomorphism of each of these groups, it suffices to show that
\begin{equation}\label{claim}
\tilde{\phi}(P) + \tilde{\phi}(\cO_k) = \tilde{\phi}(\tilde{P}),
\end{equation}
where $\tilde{P} = P + \cO_k$. 

Consider $P = (x_0,y_0)$ with $P \notin  \ker \phi$. Then we have $\tilde{P} = (x_0 + x_1k\ep, y_0 +y_1k\ep)$,  where $x_1 = -2y_0$ and $y_1 = -(3x_0^2 + A)$.  

Suppose that $t_2 \circ \tilde{\phi} = mt_1 + O(t_1^2)$. Then
\[  \begin{matrix}
\tilde{\phi}(x_0,y_0) + \tilde{\phi}(\cO_k) & = & (r(x_0), y_0s(x_0)) + (mk\ep:1:0) \\
& = & \Big(r(x_0) - 2y_0s(x_0)(mk)\ep: y_0s(x_0) - (3r(x_0)^2 + A_2)(mk)\ep: 1 \Big) 
\end{matrix}  \]

On the other hand we have
\[  \begin{matrix}
\tilde{\phi}((x_0 + x_1k\ep, y_0 +y_1k\ep)) 
& = & (r(x_0) + r'(x_0)x_1k\ep, y_0s(x_0) + (y_0s'(x_0)x_1 + y_1s(x_0))k\ep) \\
\end{matrix}  \]

Now suppose further that $P \notin E[2]$. Since the point satisfies the Weierstrass equation of $\tilde{E}_2$, the ratio of the $\ep$-coefficients of the coordinates equals $\frac{3r^2 + A_2}{2ys}$ by Lemma \ref{decomp}.  
Thus it suffices to verify that the $\ep$-coefficients of the $x$-coordinates agree. Since $x_1 = -2y_0$, this reduces to showing that $s(x_0)m = \frac{dr}{dx}(x_0)$, or equivalently, that  $\om_2 \circ \phi = m \cdot \om_1 $, where $\om_i = \frac{dx_i}{y_i} = (1 + O(t_i))dt_i$ is an invariant differential of $E_i$.  Expanding $\om_2 \circ \phi$ around $t_1$ and using the fact that $t_2 \circ \phi = mt_1 + O(t_1^2)$, we have that $\om_2 \circ \phi = (m + O(t_1))dt_1$. Thus $\frac{\om_2 \circ \phi}{\om_1} (\PI) = m$. Since $\om_2 \circ \phi$ and $\om$ are both invariant differentials, this is a constant function and (\ref{claim}) holds for $P \notin \ker \phi \cup E[2]$.

Using the equivalence of the $\ep$-coefficients of the $y$-coordinates, we have that 
\[
(3r(x)^2 + A_2)m =  (3x^2 + A_1)s(x)
\]
for infinitely many $x \in K$, and therefore this is an identity in $K(x)$. Thus, for $P \in E[2]$,  $(3r(x_0)^2 + A_2)m =  (3x_0^2 + A_1)s(x_0)$, and (\ref{claim}) holds for points of order two. 

Finally, for $P \in \ker \phi$, choose $T \notin \ker \phi$. By Definition \ref{affine} and (\ref{claim}), 
$\tilde{\phi}(P + \cO_k) =  \phi_T(P + \cO_k)  = \phi(P + \cO_k + T) + \phi(-T) = \phi(P+T) + \phi(\cO_k) + \phi(-T) = \phi(P) + \phi(\cO_k)$. Therefore, (\ref{claim}) holds for all $P \in E(K)$, and  $\tilde{\phi}$ is a homomorphism. $\Box$

\begin{lem}\label{lem2} Let $\phi: E_1 \rightarrow E_2$ be an isogeny with $t_2 \circ \phi = mt_1 + O(t_1^2)$ for $m \in K$. Then
$$e(\phi(P), \cO_{mk}) = e(P, \cO_k)^{ \deg \phi}.$$
\end{lem}

\noindent{\bf Proof:} If $\phi$ is inseparable, then the degree of inseparability is $q = p^r$ for some $r>0$ and thus $p$ divides $\deg \phi$.  Furthermore, $m = 0$ since the order of $t_2 \circ \phi$ at $\PI^{(1)}$ is 
the degree of inseparability. So both $e(P, \cO_k)^{ \deg \phi}$ and $e(\phi(P), \cO_k)^m$ equal 1, and the result holds. 

Now assume $\phi$ is separable, which implies that $m \neq 0$. By the proof of Proposition \ref{sumofslopes}, it suffices to show that
$$m\frac{df_{P_2}/dt_2}{f_{P_2}}(\PI^{(2)}) = (\deg \phi) \frac{df_{P_1}/dt_1}{f_{P_1}}(\PI^{(1)}).$$

Let $\ker \phi = \{R_1, ... R_s\}$. Since $\dv(f_{P_2}) = p(P_2) - p(\PI)$ and $\phi$ is separable, $\dv(f_{P_2} \circ \phi) = 
\sum_{i=1}^s p(P_1 + R_i) -p(R_i)$. Let $g_i = \frac{ l_{-P_1,-R_i}}{v_{P_1}v_{R_i}}$. Then $\dv(f_{P_2}\circ \phi)   = \sum_{i=1}^s p[(P_1) - (\PI) + \dv(g_i)] = \sum_{i=1}^s \dv(f_{P_1}g_i^p)$. 
Thus, up to a constant, $f_{P_2} \circ \phi = f_{P_1}^{\deg \phi} (\prod_{i=1}^s g_i)^p$. Since the characteristic of $K$ is $p$, 
$$d(f_{P_2}\circ \phi) = (\deg \phi)f_{P_1}^{\deg \phi -1}(df_{P_1})(\prod_{i=1}^s g_i)^p.$$ 
Thus
\begin{equation}\label{phifact}
\frac{d(f_{P_2} \circ \phi) }{f_{P_2} \circ \phi} = (\deg \phi)\frac{df_{P_1}}{f_{P_1}}.
\end{equation}

Note that for any function $g$ expanded around $t$,  $\frac{dg}{dt} \circ \phi = \frac{d(g \circ \phi)}{d(t \circ \phi)}$. 
Using this and  (\ref{phifact}), we have
\[  \begin{matrix}
m\frac{df_{P_2}/dt_2}{f_{P_2}}(\PI^{(2)}) & = & m \Big( \frac{df_{P_2}/dt_2}{f_{P_2}} \circ \phi \Big) (\PI^{(1)}) \\
& = & m\frac{d(f_{P_2} \circ \phi) /d(t_2 \circ \phi)}{f_{P_2} \circ \phi}(\PI^{(1)}) \\
& = & m (\deg \phi)\frac{df_{P_1}/d(t_2 \circ \phi)}{f_{P_1}} (\PI^{(1)}).
\end{matrix}  \]

\noindent From that $\frac{dt_1}{d(t_2 \circ \phi)} = m^{-1} + O(t_1)$, we have
\[  \begin{matrix}
m\frac{df_{P_2}/dt_2}{f_{P_2}}(\PI^{(2)})  
& = & (\deg \phi)  \frac{df_{P_1}/dt_1}{f_{P_1}}(\PI^{(1)}),
\end{matrix}  \]
and the lemma is proved. $\Box$ \\

The proof of Proposition \ref{isogprop} is now immediate. From Lemma \ref{lem2} and Definition \ref{isogdef}, we have 
\[  \begin{matrix}
e_p(\tilde{\phi}(P),\tilde{\phi}(\cO_k)) & = & e(\phi(P), \cO_{mk}) 
& = & e_p(P, \cO_k)^{\deg \phi}. 
\end{matrix}  \]
Thus, since $e_p$ is bilinear and $\tilde{\phi}$ is a homomorphism, the proposition holds. 

\section{Another application of elliptic curves over the dual numbers}\label{smart}

We have seen how the extension of the Weil pairing to $p$-torsion over the dual numbers directly leads to the previously defined maps of \cite{ruck} and \cite{sem}. Though we have not gained any ``new" information, we have shown that discrete logarithm attacks on $p$-torsion subgroups of \cite{ruck} and \cite{sem} may be interpreted as Weil-pairing-based attacks, exactly the same as the MOV attack on prime-to-$p$ torsion subgroups. In this section, we give another example of how looking at elliptic curves over the dual numbers may be a fruitful approach.

The DLP attack on anomalous curves of Smart \cite{smart} involves working in  $\tilde{E}(\zz/p^2\zz)$ where $\tilde{E}$ is a non-canonical lift of $E$ (meaning $p$-torsion points of $E$ are no longer $p$-torsion  when lifted to $\tilde{E}$). The attack involves lifting points $P,Q \in E[p]$ with $Q = nP$ to $\tilde{E}(\zz/p^2\zz)$, 
multiplying the points by $p$, and applying the map $(x,y) \mapsto \frac{x}{y}$. In this way, solving for $n$ such that $nP = Q$ reduces to solving an instance of the DLP in $\ff_p^+$. The fact that this map is a homomorphism may be shown via the $p$-adic elliptic logarithm (see \cite{smart}, or \cite{lcw}, p. 190).

If we consider $\tilde{E}(\dup)$ instead, the attack works analogously, and the reasoning behind it is elementary. (In fact, the attack works for $\tilde{E}(\du)$, where $K$ is any field of characteristic $p\neq 0,2,3$.) Lift $P,Q$ to  $\tilde{P}, \tilde{Q} \in \tilde{E}(\dup)$. The points $\tilde{P}, \tilde{Q}$ may no longer be dependent. However, since $nP=Q \in E(\ff_p)$, there exists $R \in \Theta$ such that $n\tilde{P}-\tilde{Q} = R$. Since $P,Q$ are points of $p$-torsion, $p\tilde{P}, p\tilde{Q} \in \Theta$. 
Thus we have the following equation in $\Theta$
\begin{equation}\label{dlp}
p(n\tilde{P}) - p\tilde{Q} = pR = \PI.
\end{equation}
Note that $p\tilde{P}, p\tilde{Q} = \PI$ if and only if $\tilde{P}$ and $\tilde{Q}$ are $p$-torsion points in $\tilde{E}$. Thus if this is not the case,
we can translate (\ref{dlp}) to an instance of the DLP in $\ff_p^+$ via the homomorphism $(k\ep:1:0) \mapsto k$ and then solve for $n$. 

This version is more efficient, as computations in $\ff_p[\ep]$ are more straightforward than in $\zz/p^2\zz$. It may present another advantage as well, related to the fact that the DLP attack requires that the lift of the curve over $\ff_p$ be non-canonical. 

Let $\tilde{E}$ be any lift of the curve $E: y^2 = x^3 + Ax + B$, with $j$-invariant $j \in \ff_p$. Note that $D = 4A^3 + 27B^2 \neq 0$ since $E$ is non-singular. Define $j(\tilde{E})$  as the value $\frac{4\tilde{A}^3 }{4\tilde{A}^3 + 27\tilde{B}^2}$. Since $D \neq 0$, the denominator is invertible, and hence $j(\tilde{E}) \in \dup$. Let $\tilde{j}$ denote the value $j(\tilde{E})$, and note that $\tilde{j} \equiv j \mod \ep$.  The following proposition shows that $\tilde{j} \in \ff_p$ if and only if the elliptic curve $\tilde{E}$ can be transformed to the ``canonical lift" (as defined in Section \ref{prelim}) by an invertible change of coordinates.  

 \begin{prop}\label{coc}
Let $E$ be given by $y^2 = x^3 + Ax + B$.  Let $\tilde{E}$ be a lift of $E$ to $\dup$ with $\tilde{A} = A + A_1\ep, \tilde{B} = B + B_1\ep$, for  $A_1,B_1 \in \ff_p$. 
Then $\tilde{j} \in \ff_p$ if and only if there exists $\mu = 1 + kt$ with $k \in \ff_p$ such that $\mu^4A = \tilde{A}$ and $\mu^6B = \tilde{B}$. In this case, there exists a change of coordinates $x \mapsto \mu^2x, y \mapsto \mu^3y$ taking $E$ to $\tilde{E}$, where $E$ is viewed as an elliptic curve over $\dup$.
\end{prop}

\noindent{\bf Proof:} Assume there exists  $\mu = 1 + kt$, $k \in \ff_p$ with $\mu^4A = \tilde{A}$ and $\mu^6B = \tilde{B}$. Then $\tilde{j} = \frac{4\tilde{A}^3 }{4\tilde{A}^3 + 27\tilde{B}^2} = j$. 

For the other implication, assume $\tilde{j}  \in \ff_p$. Then 
 $\tilde{j} = j$ and a calculation with the $\ep$-components yields
\begin{equation}\label{eq}
12A^2A_1D = 4A^3(12A^2A_1 + 56BB_1).
\end{equation}
To find $\mu = 1 + kt$ such that $\mu^4A = \tilde{A}$ and $\mu^6B = \tilde{B}$, we solve $4kA = A_1$ and $6kB = B_1$ simultaneously for $k$. If either $A$ or $B$ is zero this is no problem. If $A, B \neq 0$, choose $k \in \ff_p$ such that $4kA = A_1$. Then (\ref{eq}) becomes 
$$12A^2(4kA)D = 4A^3(12A^2(4kA) + 56BB_1)$$
which simplifies to $6k(D-4A^3) = 27BB_1$. This implies that $B_1 = 6kB$, as desired. $\Box$ \\

Thus if $\tilde{j} \in \ff_p$, the $p$-torsion of $E$ lifts to $p$-torsion of $\tilde{E}$, and the DLP attack over the dual numbers fails.  Calculations suggest that lifts with $\tilde{j} \in \ff_p$ are the only lifts of $E$ for which $p$-torsion lifts to $p$-torsion. Presuming this,  it is easy to avoid a lift to $\dup$ for which $\tilde{P}$ and $\tilde{Q}$ are $p$-torsion simply by choosing a lift with $j$-invariant $\tilde{j} \notin \ff_p$. This differs from the case of lifting to $\zz/p^2\zz$, since (to the author's knowledge) there is no analogously simple way to determine from the $j$-invariant $\tilde{j} \in \zz/p^2\zz$ whether or not the lift is canonical.

\end{document}